\documentclass[12pt]{article}
\usepackage{amscd,amsfonts,amsmath,amssymb,amstext,amsthm}

\title{\Large{\textbf {Hans Grauert (1930-2011)}}
\footnote{
AMS Subject classification 32C15, 32C55, 32E40, 32F10, 32J25, 32L99, 32Q15,
32Q15,01A61}}
\author{Alan Huckleberry}
\theoremstyle{plain}

\newtheorem*{theorem} {Theorem}
\begin{document}
\maketitle
\noindent
Hans Grauert passed away at the age of 81 in September of 2011.
His contributions to mathematics have and will be used with great frequency, 
and in particular for this reason will not be forgotten. 
All of us in mathematics stand on the shoulders of giants.
For those of us who work in and around the area of 
\emph{complex geometry} one of the greatest giants of the second
half of the 20th century is Hans Grauert. 

\bigskip\noindent
Specialists in the area know this, but even for them his
collected works, annotated with the much appreciated help of
Yum-Tong Siu, should at least be kept on the bedside table. 
An eloquent firsthand account of the \emph{Sturm und Drang}
period in M\"unster can be found in Remmert's
talk  (an English translation appears in \cite{R}) on the
occasion of Grauert receiving the von Staudt Preis
in Erlangen. More recently, on the occassion of his 
receiving the Cantor Medallion, we presented a sketch of the 
man and his mathematics (see \cite {H1,H2}). 
In the AMS-memorial article \cite{AMS} specialists
in the area, some of whom were students of Grauert, give us a closer
look. In the present article we attempt to give an in-depth view, written for
non-specialists, of Grauert's life in mathematics and the 
remarkable mathematics he contributed.
\subsection* {Early surroundings}
Grauert was born in 1930 in Haren, a small town near the Netherlands in the
northwestern part of Germany.  Many of our
friends who lived as children in this region recall their wartime fears,  
in particular of the bombings.   
M\"unster, which, together with the neighboring city Osnabr\"uck,
was the city of the signing of the 
Treaty of Westfalia ending the Thirty Years' War, 
was to a very large extent flattened.  We never
heard Grauert mention any of this; instead, he often told stories
about having fun playing with unspent shells after the war, 
something that took the sight of another great complex analyst of 
the same generation, Anatole Vitushkin, far away in the Soviet Union. 
In \cite{H1,H2} we recalled Grauert's
detailed remarks about these days at his retirement dinner.  In particular,
he wanted to explicitly thank one of his grade school teachers for not
failing him for his lack of skill in computing with numbers, informing him
that soon he would be thinking in symbols and more abstractly.

\bigskip\noindent
Immediately after completing elementary school and his gymnasium education
in nearby Meppen, Grauert began his studies in Sommersemester 1949 
at the university in Mainz. In the fall of that year he transferred back
to M\"unster, where he would go from schoolboy to one of the worldwide
leading authorities in the area of several complex variables and
holder of the Gauss Chair in G\"ottingen in a period of ten years.

\bigskip\noindent
Despite the destruction caused by the war (Germany was only beginning
to \emph{rise from the ashes}), M\"unster was one of the
best places in the world to start out in complex analysis.  At the
leadership level Heinrich Behnke, Henri Cartan and Karl Stein were
playing key roles. Among the students there were already the likes of
Friedrich Hirzebruch, who began his studies in 1945, and Reinhold Remmert,
who would become Grauert's lifelong friend and co-author of numerous
fundamental research articles and expository monographs.
  
\bigskip\noindent
Behnke had come to M\"unster in 1927 as a proven specialist
in the complex analysis in several variables of the time. 
Fundamental first results had already been discovered and proved.
These include the remarkable facts about the location and nature
of singularities of holomorphic functions proved by two giants of 
the early 20th Century, Eugenio Elie Levi and Friedrich Hartogs.
To give a flavor of the times, let us summarize a bit of this
mathematics. 

\bigskip\noindent
Levi had understood that if the smooth boundary
of a domain $D$ in $\mathbb {C}^n$, $n\ge 2$, is locally defined
as $\{\rho=0\}$ with $\rho $ being negative in the domain, then
the complex Hessian $(\frac{\partial^2\rho}{\partial z_i\partial \bar{z_j}})$
contains curvature information which determines whether or not a holomorphic
function on the domain can be continued holomorphically across the
boundary point in question.  For example, the appropriate notion
of concavity at a boundary
point $p$ with respect to this \emph{Levi form}
can be restated by requiring the existence of a holomorphic
mapping $F$ from the unit disk $\Delta $ in the complex plane to the closure
of $D$ with $F(0)=p$ and $F(\Delta \setminus \{p\})\subset D$. If 
$\partial D$ is concave at $p$, then every function holomomorphic 
on $D$ extends holomorphically to a larger domain $\hat {D}$ which
contains $p$ in its interior. 

\bigskip\noindent
Hartogs had understood related phenomena,
the simplest example of which goes as follows. Consider a domain
$D$, e.g., in $\mathbb {C}^2$ which can be viewed as a fiber space by the
projection onto the unit disk $\Delta $ in the first variable.
For some arbitrarily small neighborhood $\Delta'$ of the origin
in $\Delta $ the fibers are assumed to be unit disks in the
space of the second variable and otherwise are annuli with
outer radius $1$ and inner radius arbitrarily near $1$.  Then
every holomorphic function on $D$ extends holomorphically to
the full polydisk $\hat {D}=\{(z_1,z_2):\vert z_i\vert<1, \ i=1,2\}$.

\bigskip\noindent
The so-called Cousin problems, formulated by Cousin in a special context
in the late 19th century, which when positively answered
are the higher-dimensional analogues of the Mittag-Leffler 
and Weierstrass product theorems of one complex variable,
are also related to questions of analytic continuation.  For example,
Cousin I for a domain $D$ in $\mathbb {C}^n$, asks for the
existence of a globally defined meromorphic function on $D$ with
(locally) prescribed principal parts. This means that on an open covering
$\{U_\alpha\}$ of $D$ there are given meromorphic functions $m_\alpha $
which are compatible in the sense that $m_\beta -m_\alpha =:f_{\alpha\beta}$
is holomorphic on the intersection $U_{\alpha\beta}$.  The question
is whether or not there is a globally defined meromorphic function
$m$ on $D$ with $m-m_\alpha$ holomorphic on $U_\alpha$ for every $\alpha$.

\bigskip\noindent
The following is a connection of Cousin I to the study of 
analytic continuation. Let $D$ be a domain in $\mathbb C^n$ 
with $p=0$ in its smooth boundary and suppose that the set
$\{z_1=0\}$ is locally in the complement of $D$.  In order to show
that some holomorphic function on $D$ cannot be continued across
$p$, one could try the following: 
Bump out $D$ at $0$ to obtain a slightly larger
domain $\hat D$ which contains an open neighborhood $U_0$ of $p$.
Define $U_1=D$ and consider the Cousin I data for the covering
$\{U_0,U_1\}$ of $\hat D$ of $m_1\equiv 0$ and $m_0=\frac{1}{z_1}$.
If this has a ``solution'' $m$, then $m\vert D$ is an example of
a holomorphic function that cannot be continued through $p$!  We
mention this here, because, as we explain
below, Grauert brought this Ansatz to fruition and perfected it 
in its ultimate beauty.

\bigskip\noindent
Behnke certainly knew that several complex variables was an
area ripe for development and set about building a research
group for doing so.  He was an active mathematician who understood where
mathematics was going and where it should go. He was optimally
connected to the world outside M\"unster. Caratheodory,
Hopf, Severi and many others were his close friends. Perhaps
above all, he was a remarkable organizer of all sides of
our science!  It must be emphasized, however, that he was fortunate!
Even early on he had a group of magnificent students/assistants, three of
whom we have had the honor of knowing: Peter Thullen, Friedrich Sommer
and Karl Stein. Secondly, for the seemingly innocuous reason that
he had proved a small remark on circular domains which improved on
an old result of Behnke, in 1931 Henri Cartan was invited by
Behnke to give a few talks in M\"unster. 

\bigskip\noindent
Thullen and Cartan became good friends and proved a basic result 
characterizing \emph{domains of holomorphy} which possess holomorphic functions
which cannot be continued across their boundaries. Thullen went
on to prove a number of results, including important continuation 
theorems. Behnke and Thullen published their Ergebnisbericht which
in particular outlines the key open problems of the time. In
a series of papers which are essential for certain
of Grauert's works, Kyoshi Oka solved many of these problems.  Story has
it that he wrote Behnke and Thullen a thank-you note for posing
such interesting questions.

\bigskip\noindent
Friedrich Sommer, who was one of the founding fathers of the
Ruhr University and who was responsible for continuing the tradition
of complex analysis in Bochum, was one of the stalwarts of
the Behnke group which Grauert joined.  
 
\bigskip\noindent
Before the war, Behnke was still active in mathematics research,
in particular with Stein who after the war 
became the mathematics guru of complex
analysis in M\"unster. (See \cite{H3} for a detailed discussion of
Stein's contributions.)  The work of Behnke and Stein underlining
approximation theorems of Runge type, and, for example,
Stein's emphasis on implementing concepts from algebraic topology
(he spent time in Heidelberg with Seifert) certainly influenced
the young Grauert.

\bigskip\noindent
Not enough can be said about the importance of Henri Cartan
for the \emph{M\"unsteraner school of complex analysis}. The pre-war 
interaction indicated above was just the beginning.  
Despite the fact that the Nazi atrocities 
directly touched Cartan's family (his brother was 
assasinated in 1943), shortly after the war, in 1947, he accepted
Behnke's invitation to visit M\"unster. For those of us who did 
not experience the horrific events of that time, it is difficult to
imagine the magnitude of importance, maybe most importantly 
at the human level,
of Cartan's reestablishing the Paris--M\"unster connection (see \cite{HP}
for more on the importance of Cartan for postwar German mathematics). The
importance for complex analysis, in particular for Hans Grauert,
is discussed below.
       
\subsection* {Initial conditions}
When Grauert arrived in M\"unster, despite the fact that the
worldly ammenities of the university were still at best minimal,
Behnke had complex analysis up and running and, in a certain sense,
the conditions for research were optimal. On the one hand there
was Stein, a kind, modest man of enormous enthusiam and 
energy who had deep insight
at the foundational level of, e.g., analytic sets, holomorphic mappings, etc.
Certainly the Cousin problems and their relationships to domains
of holomorphy had guided a big part of his thinking. As a result of
research with Behnke before the war and published in 1948, he knew
these were solvable on non-compact Riemann surfaces. To add a bit
to a paper which he worried was otherwise too short he formulated
three axioms for what are now known as \emph{Stein manifolds}
which he felt would be the correct general context for solving problems
of Cousin type (\cite{S}): Globally defined holomorphic functions
separate points and give local coordinates, and given a divergent
sequence $\{z_n\}$ there exists a holomorphic function $f$
with $\vert f(z_n)\vert$ unbounded. 

\bigskip\noindent
Stein was a hands-on craftsman and this certainly
influenced the spirit of the Behnke seminar where there were
lengthy naive (healthy!) discussions of examples such as $\sqrt{xy}$ (the cone
singularity $z^2-xy=0$ which can be viewed as a 2:1 cover of
$\mathbb {C}^2$ ramified only at the origin). On the other hand,
the mathematics world outside M\"unster, particularly influenced
by developments in France, had made a quantum leap in sophistication.
However, Behnke made sure that M\"unster was not isolated.

\bigskip\noindent
Fritz Hirzebruch had begun his studies in M\"unster in 1945.
He lived in Hamm where Stein also lived.  We have heard that
they traveled to M\"unster by train together often hanging on
to the outer running boards with Stein making propaganda for
the role of algebraic topology in complex geometry.  During
Hirzebruch's studies, Behnke sent him to his friend Hopf
in Zurich.  Hirzebruch happily reminisced about learning
from Hopf about blowing up points and blowing down curves in surfaces.
In fact his thesis (published in 1951) is a jewel about
surfaces where this process plays a role. On another not unrelated
topic, in a talk to a historical society on the Riemann-Roch
theorem, Hirzebruch said that probably the most important new development
for him in the early 1950's was understanding of the notion of a line bundle!
Just a few years later he fused a hefty portion of the new 
sophistication with his own ideas to prove his Riemann-Roch Theorem!!
(published in 1954.) For the young M\"unsteraner it
must have been extremely motivating to see this remarkable development.

\bigskip\noindent
Cartan's early works, e.g., with Thullen and those on automorphism groups
of domains, fit in the style of complex analysis at
the time. However, Cartan not only made the leap from the classical 
to the post-war level of sophistication, he was one of the main 
figures who shaped it. Despite having formulated and proved
Theorems A and B (published in 1951), which not only solve 
the Cousin problems on Stein spaces but put complex analysis 
in another world of abstraction (the distance from M\"unster to 
Paris could no longer be measured in kilometers), he remained
in contact with and supported the members of the Behnke group.
By the way, it was his idea to refer to these spaces as
\emph{Vari\'{e}t\'{e} de Stein}.

\subsection* {Crescendo}
At this point in the historical timeline Grauert entered the picture
and, at certain points with the help of distinguished co-workers, took
complex analysis to yet another level. Having set the stage above, we
now turn to a description of representative aspects of his published
works. We begin with an overview.

\bigskip\noindent
Grauert received his doctorate in M\"unster in mid-1954.  His first
publications appeared in 1955, the publication from his thesis in
1956. In the five or six years that followed, his contributions
to mathematics were truly remarkable: a wealth of ideas,
numerous basic results and simply quite a number of published
pages. Disregarding research announcements (Comptes Rendus Notes),
conference reports and expository articles, in this intense period 
he authored or co-authored (with Remmert, Andreotti and one with
his student Docquier) 19 articles which covered a total of roughly 600
pages. In the three or four years after this period, when both
he and Remmert were in G\"ottingen, they jointly wrote three basic 
research monographs in book form: \emph{Analytische Stellenalgebren}, 
\emph{Stein Theory} and \emph{Coherent
Analytic Sheaves}. The final versions of the latter two appeared
much later.  After settling in G\"ottingen, where he also devoted
a great deal of time to his students (he guided more than 40 Ph.D.
theses), Grauert continued to make important research contributions.   
Altogether he published more than 90 works, most
of which were devoted to topics in the areas of several complex variables
and complex algebraic geometry.

\bigskip\noindent
In a nutshell one can summarize Grauert's work as being fundamental
for the foundations of the geometric side of complex analysis,
particularly his early work with Remmert, and for our understanding of
the multifaceted global phenomena related to Levi curvature. 
His solution of a certain Levi problem 
is just one of a number of results in this direction. There are two 
early works of Grauert that stand out as the peaks among many
mountains: The Oka Principle (1957) and the Direct Image Theorem (1960). 
These and selected works in the areas indicated above are 
discussed in some detail in the next section.

\bigskip\noindent
In addition to those works which will be discussed 
in the next section, a number of important papers must be 
mentioned, e.g., that on the solution of the
Mordell conjecture in the function field case. 
Weil mused that Manin, the algebraist, used analytic methods for 
this whereas Grauert, the analyst, approached it algebraically. 
In fact, if one looks at the
paper, one immediately sees Grauert's geometric viewpoint. Other
results which stand out are his construction of
the versal deformation space for compact complex spaces (simultaneously
with Douady) and that for deformations of isolated singularities, 
his basic cohomology vanishing theorem with Riemenschneider and results 
on conditions for the formal equivalence of neighborhoods 
of analytic subsets implying convergent equivalence. His work with
M\"ulich on vector bundles on $\mathbb {P}_2$ has been extremely 
influential. Fundamental work on the analytic side, in particular solving
$\bar {\partial}$-problems with bounded data, was carried out
with his students, Ramirez and Lieb.  He also wrote 
textbooks for basic real analysis with Lieb and for
linear algebra with Grunau.  In the area of several complex variables
he wrote two textbooks with Fritzsche, together with Peternell and Remmert
he edited and contributed several chapters to a volume of the 
Encyclopedia of Mathematical Sciences and wrote the three research monographs
with Remmert which were mentioned above.

\bigskip\noindent
Grauert considered a wide range of topics.  For example, one should 
not forget his ideas on hyperbolicity (see Demailly's comments in \cite{AMS}) 
as well as his interests in vector bundles, deformation theory and 
in understanding analytic equivalence relations, a 
topic that had followed him since his early encounters with Karl Stein.  
He had a philosophical side as well which went along with his
desire to understand certain kinds of physical (quantum mechanical) 
phenomena. It seems that he read Riemann's work having this in mind and, 
based on this, developed his own theory of discrete geometry. We recall 
his series of lectures at Notre Dame on his axiomatic approach and 
note that at the end of Volume II of his collected works he 
included several pages on this.
Given that he obviously carefully polished these two volumes, it
is clear that he took this subject very seriously and that it meant
a great deal to him.

\subsection* {Comments on selected works}
Under the given time and space constraints it is only possible
to present a small sample of Grauert's works.  Since he is perhaps
best known for his results at the foundational level in complex geometry,
those involving Levi geometry, his Oka principle
and his proof of the direct image theorem, our remarks here will 
focus on these subjects.
\subsubsection* {Early days}
We begin with comments on certain aspects of Grauert's
dissertation (published in \cite{GDiss})
which underlined an important connection between
complex differential geometry and complex analysis.  This is followed by
remarks on the paper where he vastly improved our
foundational understanding of Stein spaces (\cite {GStein}).  
This early work emphasized the need for building the foundations 
of \emph{complex spaces}.  Much in this direction was accomplished in the
basic paper \emph{Komplexe R\"aume} (\cite {GR1}) of Grauert and 
Remmert which is the third paper we review in this paragraph.

\bigskip\noindent
\textbf{Charakterisierung der Holomorphiegebiete durch die
vollst\"andige K\"ahlersche Metrik}

\medskip\noindent
When Grauert  went to the ETH (1953), it was already quite fashionable to
study K\"ahler manifolds.  By definition such a manifold possesses a
Hermitian metric whose imaginary part is a closed 2-form.
Locally this form $\omega $ has a potential $\varphi $, i.e.,
$\omega =\frac{i}{2}\partial\bar{\partial}\varphi $, and the
positivity of the metric translates to $\varphi $ being strictly
plurisubharmonic. Hodge's book, which appeared in 1941, was
well known and Eckmann and Guggenheimer were busy in Zurich (Guggenheimer
went to Israel in 1954) looking at more general manifolds. The fact
that plurisubharmonic functions (Lelong, Oka 1942) are important in
complex analysis was widely understood.  K\"ahler himself had 
thought in terms of the potential function and had in fact 
proved that K\"ahler is, as 
mentioned above, the
same thing as having a locally defined strictly
plurisubharmonic potential. Much was in the air when Grauert 
started thinking in this direction.

\bigskip\noindent
At the beginning of his paper Grauert states that it is ``naheliegend'' to
study the connection between complete K\"ahler metrics and domains
of holomorphy $D$. In hindsight this is true, but he was the
first to do this. For our purposes a \emph{domain of holomorphy}
is a domain in $\mathbb {C}^n$ for which there exists a function
$f$ holomorphic on $D$ which cannot be extended to a function
holomorphic on a larger manifold. In fact, given a sequence $\{z_n\}$
which converges to every boundary point, one can construct $f$
with the property that $\lim \vert f(z_n)\vert=\infty$. 

\bigskip\noindent
Grauert begins the article by pointing out that, given a complete
K\"ahler metric on $D$ and a closed analytic subset $A$ of $D$, i.e.,
a set defined as the common $0$-set of finitely many holomorphic
functions, it is a simple matter to adjust the given K\"ahler
metric appropriately to obtain a complete K\"ahler metric on
$D\setminus A$.  Since holomorphic functions extend across analytic
sets which have codimension at least two, it is then immediate that
there are domains with complete K\"ahler metrics which are not
domains of holomorphy.

\bigskip\noindent
After making the above remark, Grauert then proves that the desired
result holds if $D$ has a $\mathbb {R}$-analytic boundary.  In other words,
such domains are domains of holomorphy if and only if
they possess a complete K\"ahler metric. This result, which required
substantial technical work, stimulated
a great deal of futher research, first of all due to the idea
that the existence of a complete K\"ahler metric and pseudoconvexity
are related. The regularity question also turned out to be of interest.  
For example, some years later Ohsawa showed that only a $C^1$-boundary 
is necessary (\cite{O}).

\bigskip\noindent
This paper is Grauert's doctoral thesis.  He profusely thanks Behnke 
and Eckmann. I would guess that Eckmann and Guggenheimer 
discussed K\"ahlerian geometry with him, although their work exclusively dealt
with compact manifolds and in particular had nothing to do with
pseudoconvexity.  Grauert received his
degree in M\"unster in July of 1955, the first referee being Behnke and the
second Friedrich Sommer.

\bigskip\noindent
\textbf{Characterisierung der holomorph vollst\"andigen
komplexen R\"aume}

\medskip\noindent
In this paper Grauert is clearly fascinated by the question of 
countability of the topology
of complex spaces (Rado noticed this in the case of Riemann surfaces). 
Throughout the paper a complex space is an $\alpha $-space and Grauert
is thinking in terms of it locally being the graph of a multivalued 
holomorphic function.  He does not yet have
the result that every $\alpha $-space is a $\beta_n$ space
(see our review of \emph{Komplexe R\"aume} for the notation). As a result
he proves his main results under assumption $C$ which due to the later
work of Grauert and Remmert just means that a complex space is locally
the common $0$-set of finitely many holomorphic functions on a domain
in $\mathbb {C}^n$.  

\bigskip\noindent
The \emph{new} axiom here is that of $K$-Vollst\"andigkeit, i.e., that
global holomorphic functions define a map at a given point 
which is finite fibered
near the point in question.  Much work then shows that under this
condition the topology is countable and finally that the 
$n$-dimensional space $X$ is globally a ramified Riemann domain 
over $\mathbb {C}^n$, i.e., that there is a generically maximal rank
holmorphic map $F:X\to \mathbb{C}^n$.  Then, using nontrivial 
(but more or less classical) methods, Grauert shows that if
$X$ is holomorphically convex, then it is Stein. 
The main result then is an essential weakening of Stein's axioms:
$K$-vollst\"andig plus holomorphic convexity are equivalent to the following
four axioms of Stein:
\begin {enumerate}
\item
Countable topology.
\item
Globally defined
holomorphic functions give local embeddings.
\item
Globally defined
holomorphic functions separate points. 
\item
Holomorphic convexity,
i.e., given a divergent sequence $\{x_n\}$ there exists a holomorphic
function $f$ with $\lim \vert f(z_n)\vert =\infty$.
\end {enumerate}

\bigskip\noindent
\textbf{Komplexe R\"aume}

\bigskip\noindent
Here the authors work from the point of view of the definition
of Behnke and Stein which is that a complex space $X$ is locally 
the graph of a multivalued holomorphic function on $\mathbb {C}^n$. 
This was nothing new in the 1-dimensional case, because under the
assumptions of Behnke and Stein the resulting space is smooth and
locally just the graph of an algebraic function!  However, in
the higher dimensional case singularities arise.  Members of 
the Behnke seminar in those years have told us that they spent
great energy trying to understand $\pm\sqrt{xy}$ which is the
cone defined by $z^2-xy=0$ over the $xy$-plane!

\bigskip\noindent
To be precise, a Behnke-Stein complex space is a Hausdorff space $X$
which satisfies the following local condition: Every $x\in X$ is 
contained in an open neighborhood $U$ which is equipped with
a continous map $\varphi :U\to V$ onto an open set in $\mathbb {C}^n$
which contains a proper analytic subset $A$ with the property that
the restriction $\varphi :U\setminus \varphi^{-1}(A)\to V\setminus A$
is a proper finite covering map. In particular, $\varphi $ is a local
homoeomorphism and gives local holomorphic coordinates
on $U\setminus \varphi^{-1}(A)$.  The holomorphic functions on $U$ 
are then defined to be the continuous functions which are holomorphic
on $U\setminus \varphi^{-1}(A)$.  One should mention in this context
that one of the most quoted theorems of Grauert and Remmert is
that if $Y$ is, for example, a complex manifold which contains a
proper analytic subset $A$ and $F:X\to Y\setminus A$ is a proper
unramified holomorphic map, then $X$ can be (uniquely) realized
as the complement of a proper analytic subset $B$ in a larger complex
space $\bar {X}$ so that the analytic cover $X\to Y\setminus {A}$ 
can be extended to a proper (finite) holomorphic map 
$\bar {X}\to Y$ with $\bar {X}\setminus B\to Y\setminus A$ being
the original map.

\bigskip\noindent
The main goal of this paper is to show that the definition of 
Behnke and Stein is equivalent to that of Cartan and Serre of
a \emph{normal} complex space. To ``clarify'' matters Grauert 
and Remmert introduce a rather cumbersome notation.  First, the
Behnke-Stein spaces are called $\alpha$-spaces.  The spaces
coming from Paris are called $\beta $- and $\beta_n$-spaces.
The former is locally the common $0$-set $A$ of finitely many
functions on a domain $D$ in some $\mathbb {C}^n$ (of course depending
on the point) and the sheaf of holomorphic functions is just
the quotient $\mathcal {O}_D/\mathcal {I}_A$ of the sheaf of germs
of holomorphic functions on $D$ by the full ideal sheaf of functions
which vanish on $A$.  In modern terminology these are just called
reduced complex spaces.  The $\beta_n$ spaces are those which
are \emph{normal} in the sense that if a meromorphic germ satisfies
a monic polynomial equation with holomorphic coefficients, then it is itself
holomorphic.  Due to applications
it was and still is important to understand the relations among these
concepts. The Grauert-Remmert paper clears this up completely.  
Furthermore several basic results of independent interest are proved.

\bigskip\noindent
It is not terribly difficult to show that if the ramification on
an $\alpha $-space is given by a multivalued function (algebroid condition), 
then that space is a $\beta_n$-space.  That, then, is the main theorem
of the paper: $\alpha $ spaces are automatically algebroid.  Since
$\beta_n$-spaces are easily seen to be $\alpha $-spaces, this now
completes the circle: The Behnke-Stein spaces are exactly the normal
complex spaces defined by Cartan!

\bigskip\noindent
Two basic results which we have not yet mentioned were proved 
along the way: 1.) A $\beta $-space is normal if and only if the
Riemann extension theorem holds.  2.) The normalization of a 
$\beta $-space is constructed. By \emph{Riemann extension} we mean
that if $A$ is a proper analytic subset of $X$ and $f$ is holomorphic
on $X\setminus A$ and is locally bounded near $A$, then it extends
to a holomorphic function on $X$.  The normalization 
$\pi :\hat {X}\to X$ of a complex $\beta $-space, is a finite
(proper, surjective) holomorphic map 
from a canonically determined normal complex 
space which is biholomorphic at least outside of the singular set
of $X$. For $x\in X$ the number of points in $\pi^{-1}(x)$ is the
number of local irreducible components of $X$ at $p$.  It is quite
possible that the normalization is a homeomorphism with the
only difference between $\hat {X}$ being that the structure on $\hat {X}$
is richer.   

\bigskip\noindent
The reader should consult 
\emph{Coherent Analytic Sheaves} \cite{GR2} for a more modern formulation 
of the above.  On the other hand, this book was written in note form
in G\"ottingen in the early to mid-1960s, i.e., not long after
the original article.

\subsubsection* {Oka principle}
Grauert wrote three papers (\cite{Oka1,Oka2,Oka3}) on what is now 
called Grauert's Oka principle. The first two are full of 
deep ideas and fundamental work.  The third is devoted to a 
statement of what was proved in the first two along with several 
applications.  After giving some background on what was known 
when Grauert entered the picture, 
with the help of Cartan's formulation and streamlining we outline the
statements and ingredients of proof of Grauert's results. These
can be seen as proving that on a Stein space the only obstructions
to solving problems of a complex analytic nature are topological. 

\bigskip\noindent
\textbf{Approximationss\"atze f\"ur holomorphe Funktionen mit
Werten\\ 
in komplexen R\"aumen}

\medskip\noindent
\textbf{Holomorphe Funktionen mit Werten in komplexen Lieschen\\
Gruppen}

\medskip\noindent
\textbf{Analytische Faserungen \"uber holomorph-vollst\"andigen
R\"aumen}

\bigskip\noindent
The following is a simple but fundamental example
of the Oka principle.  Let $G$ be a domain in $\mathbb {C}^n$
and $D$ be a divisor on $G$, i.e., $D$ is given locally on
a covering $\{U_i\}$ by meromorphic functions $m_i$ which satisfy
the compatibility condition that $m_i=f_{ij}m_j$ on the intersection
$U_i\cap U_j=:U_{ij}$
with the $f_{ij}$ being nowhere vanishing holomorphic functions on the 
intersections $U_{ij}=U_i\cap U_j$.  One might ask (the second
Cousin problem) if there is a globally defined meromorphic function $m$
on $D$ with this divisor.  In other words, $m$ would be required
to have the same poles and zeros (counting multiplicity) as the
$m_i$ in the sense that the functions $\frac{m}{m_i}=:f_i$ are 
holomorphic and nowhere vanishing on the $U_i$.  This is of course
the same as asking for the existence of such $f_i$ with
$f_im_i$ being globally defined.

\bigskip\noindent
On the complex plane, i.e., for $D=\mathbb {C}$, the question of 
existence of the global function is answered in the positive by
describing one such as the quotient of Weierstrass products.
As a consequence of the Theorem of Behnke and Stein this even
holds for every non-compact Riemann surface.  On the other hand
for compact Riemann surfaces, and of course for higher-dimensional
compact complex manifolds, such a theorem does not hold and the
fact that it does not guides many questions in the theory.

\bigskip\noindent
Returning to the non-compact case, if $D=\mathbb {C}^n$, the
second Cousin problem can also be answered in the positive,
even by using methods that are analogous to the Weierstrass products.
However, it was realized early on that without further assumptions,
even for domains $D$ in $\mathbb {C}^n$, there would be no hope of 
solving this problem.  The appropriate class of domains which
appeared natural for solving this problem is the class of Stein
domains or domains of holomorphy.  Such is the natural domain
of existence of some holomorphic function $f$ in the sense
that it cannot be extended holomorphically to any larger complex
manifold.  Although Stein domains are optimal from many points
of view, as was realized by Oka, even on Stein domains there
are obstructions to solving such problems.  The point is
that there might not even exist continuous functions $f_i$ with
this property. In our more modern language, the $f_{ij}$ define
a holomorphic line bundle $L(D)$ on $D$ and the Cousin II problem
has a positive solution if and only if $L(D)$ is holomorphically
trivial. The problem has a continuous solution if and only if
$L(D)$ is topologically trivial. Oka's basic theorem, proved
in the pre-war years, states that a holomorphic line bundle
on a Stein domain is holomorphically trivial if and only if
it is topologically trivial. Cartan's Theorem B, or the more
general theorem of Grauert (see above), immediately implies
this statement on Stein spaces.  It should be remarked that
there are no topological obstructions to solving 
the additive (Cousin I) problem, i.e., that which asks for
a globally defined meromorphic function with prescribed 
principle parts.  On a Stein space it always has a positive
solution. 

\bigskip\noindent
For the reason sketched above, and for various other questions which
arise in complex analysis, a vague Oka principle can be formulated: 
A problem on a Stein space which is formulated in
complex analytic terms has a complex analytic solution if and only
if it has a topological solution.  This was more or less the
state of the theory when Grauert entered the picture with his three
papers which were published in 1957. Although 
it does no justice to his work, one simple-to-state consequence
of Grauert's Oka principle is that on a Stein space the 
mapping which forgets complex structure defines an isomorphism
between the categories of holomorphic and topological vector bundles.
In other words, Oka's theorem holds for arbitrary rank.

\bigskip\noindent
Unlike the case of divisors where the transition matrices
$f_{ij}$ have values in an Abelian group and the usual cohomological
technology on Stein spaces can be applied, in Grauert's
non-Abelian setting classical results of the time cannot be
applied.  On the other hand it was certainly clear that 
in order to handle the Oka principle, e.g.,
some version of the classical Runge approximation theorem would be necessary.
In fact an extremely deep version
of this approximation theorem, one which involves homotopies of holomorphic
and continuous maps, would be of essential importance in Grauert's
theory.

\bigskip\noindent
Grauert opens his first paper as follows: 
\emph{In the present paper functions $F(r)$ from a complex
space $R$ with values in an arbitrary complex space $W$ 
will be studied.} Of course he has in mind Stein spaces;
so let us assume that $R$ is Stein and consider a 
Stein space $\check{R}$ which contains $R$.  The main
question is if such a ``function''  $F$ can be approximated
(uniformly on compact subsets) by a holomorphic function
from $\check{R}$ to $W$.  For usual functions with values in
$\mathbb {C}$ the Behnke-Stein theorem was available: The
approximation theorem holds if and only if $R$ is \emph{holomorph
ausdehnbar} to $\check{R}$. This is a condition which is
a bit complicated to formulate.  It is a substantially weakend
version of there being a continuous increasing family $R_t$, $0\le t\le 1$,
of Stein domains with $R_0=R$ and $R_1=\check{R}$. In any case,
since this condition is already necessary and sufficient for
usual functions, it was clear to Grauert that he should assume it
for his more general question.  This being clear, Grauert jumps
to a suitable setting.
 
\bigskip\noindent
Following his notation, a Lie group bundle over a complex space $R$
is a holomorphic fiber bundle $L^*(R,L)$ with fiber a complex 
Lie group $L$ and structure group $L^*$ contained in the group
of Lie group automorphisms of $L$. Note that the case of
$L=(\mathbb {C}^n,+)$ is that of a holomorphic vector bundle.
Observe that such a bundle has the identity section, fiberwise
multiplication makes sense, one has the associated Lie algebra bundle
where the exponential is biholomorphic near its $0$-section and
cohomology concepts for the sheaf of sections make sense using the
group multiplication. Note also that such a bundle is \emph{not}
a principal bundle.

\bigskip\noindent
Before going further we would like to simplify the notation and
follow that of Cartan's paper (\cite{C}) where he explains
Grauert's work with great elegance. He simply denotes the Lie group
bundle by $E\to X$ and then introduces the notion of an
$E$-principal bundle defined by a cocycle $\{f_{ij}\}$ acting on
$L$ on the left.  In this way one has the fiberwise action
$F\times_XE\to F$ on the right. With this in mind one of 
Grauert's main theorems can be formulated as follows.

\bigskip\noindent
Let $\mathcal {E}_c$ be the sheaf of continuous sections and
$\mathcal {E}_a$ be the sheaf of holomorphic sections of the
Lie group bundle $E$.
\begin {theorem}
The inclusion $\mathcal {E}_a\hookrightarrow \mathcal {E}_c$ induces
an isomorphism 
$$
H^1(\mathcal {E}_a)\cong H^1(\mathcal {E}_c)\,.
$$
\end {theorem}
Of course one of the main issues to be handled is that of 
a Runge theorem. Here is Grauert's Runge theorem for beginners.
\begin {theorem}
Let $R$ and $\check{R}$ be Stein spaces with $R$ holomorph
ausdehnbar to $\check{R}$ and let $E$ be a Lie group
bundle on $\check{R}$.  Then a holomorphic section of $E$
over $R$ can be approximated by holomorphic sections of
$\check{R}$ if and only if it can be approximated by continuous
sections on $\check{R}$.
\end {theorem}
This should be regarded as a mini-Runge theorem, because the 
final version must be proved in a context where homotopy
is involved.  Let us leave this in Grauert's language.
\begin {theorem} Ist $\mathfrak {R}$ ein holomorph-vollst\"andiger
Raum, so gibt es zu jeder in $\mathfrak {R}\times \mathfrak {T}_1$
definierten $(e,h)$-Funktion mit Werten in einem Faserraum 
$L^*(\mathfrak {R},L)$ eine $(e,h,c)$-homotope $(e,h^0)$-Funktion.
\end {theorem}
Here is Cartan's formulation.\\
Notation (The $(N,H,K)$-sheaf $\mathcal F$): 
\begin {itemize}
\item
$K$ is an auxilliary compact parameter space with
$N\subset H\subset K$ such that $N$ is a deformation retract of $K$.
\item
$\mathcal {F}(U)$ is the topological group of continuous sections 
$s(x,t):U\times K\to E(U)$ which are the identity section for $t\in N$ 
and holomorphic for $t\in H$.
\end {itemize}
\begin {theorem} If $X$ is Stein, then
\begin {enumerate}
\item
$H^0(X,\mathcal {F})$ is arcwise connected.
\item
If $U$ is holomorphically convex in $X$, then image
of the restriction $H^0(X,\mathcal {F})\to H^0(U,\mathcal F)$
is dense.
\item
$H^1(X,\mathcal F)=0$
\end {enumerate}
\end {theorem}
Of course we have swept a great deal of work under the table.
Implementing in particular the refined Runge theorem,
Grauert spends a great deal of time solving the non-Abelian 
Cousin problems.

\bigskip\noindent
Here are some of the consequences of Grauert's remarkable work.
\begin {itemize}
\item
Every continuous section $s:X\to F$ is homotopic to a holomorphic section.
\item
If two holomorphic sections are homotopic in the space of
continuous sections, then they are homotopic in the space of 
holomorphic sections.
\item
Every continuous isomorphism between $E$-principal bundles
$F$ and $G$ is homotopic to a holomorphic isomorphism.
\end {itemize}
Note that for the last result it is important that
the quotient $E_g:=(G\times_XG)/E$ is a Lie group bundle and
the bundle of isomorphisms from $G$ to $F$ is the $E_g$-principal
bundle $(F\times_XG)/E$.

\bigskip\noindent
There are numerous consequences of these results, and
even recently there has been a big explosion of further 
developments (see \cite{F}).  A big additional step, 
both technically and conceptually, was Gromov's $h$-principal (\cite{G}).

\subsubsection* {Levi convexity and concavity}
Grauert made numerous fundamental contributions to understanding
the role of Levi-curvature in complex analysis. If for example
a domain $D$ with smooth boundary in a complex manifold is defined by a 
smooth function, $D=\{\rho <0\}$, it has been known since the
beginning of the 20th century that signature invariants of
the complex Hessian of $\rho $ along $\partial D$ play an essential
role in determining the complex analytic nature of $D$. Grauert
solved the version of the \emph{Levi problem} (\cite {GL})  
which states that
if the appropriate curvature form is positive-definite, i.e.,
$D$ is strongly pseudoconvex, then $D$ is essentially a Stein manifold. 
We review this basic paper here along with four other works
where convexity, concavity or both are essential ingredients.

\bigskip\noindent
In his paper with Docquier (\cite{GD}) it is shown that if $D$ is contained
in a Stein manifold and is weakly pseudoconvex (in any of a variety
of ways) at the boundary, then it is Stein. (Oka and others had
shown this for domains in $\mathbb {C}^n$.) Andreotti and Grauert
wrote two very interesting papers where concavity is involved.
One can be regarded as a mixed signature version of Grauert's
previously handled positive-definite case where higher cohomology
spaces replace function spaces (\cite{AGr1}).  In the other they describe the
structure of the field of meromorphic functions on a pseudoconcave
space and show how to apply their methods to situations where the 
space at hand is a discrete group quotient, e.g., where the meromorphic 
functions arise as quotients of modular forms (\cite{AGr2}).  

\bigskip\noindent
Finally, we review Grauert's beautiful paper \emph{\"Uber Modifikationen
und exzeptionelle analytische Mengen} (\cite{GMod}).  
In brief, here he shows us how
to use strong pseudoconvexity in the theory of compact complex spaces,
in particular to settings of algebraic geometric interest.  The title
indicates one of the themes in the paper where he proves that a 
compact complex subvariety in a complex space can be blown down
if and only if its normal bundle satisfies a natural curvature condition.
This is just one of many other results which are proved, e.g.,
new ampleness criteria, Kodaira type embedding theorem, etc..    

\bigskip\noindent
\textbf{On Levi's Problem and the Imbedding of 
Real-Analytic Manifolds}

\medskip\noindent
Let us begin here with a classical situation where $D$ is
a domain with smooth boundary in $\mathbb {C}^n$.  Every
point $p\in \partial D$ has an open neighborhood $U$ which
is equipped with a smooth function $\rho $ with nowhere vanishing
differential so that $U\cap D=\{\rho <0\}$.  In particular,
$\partial D\cap U=M$ is the smooth hypersurface $\{\rho=0\}$.  
Note that the full (real) tangent space $T_pM$ contains a
unique maximal complex subspace $T_p^{CR}M$, the \emph{Cauchy-Riemann
tangent space} to $M$ at $p$, which is $1$-codimensional over $\mathbb {R}$.
Let us say that the Levi-form $L_p(\rho)$ of the defining function 
$\rho $ at $p$ is the restriction
to the CR-tangent space of the complex Hessian
$$
Hess_p(\rho)=\Big(\frac{\partial^2}{\partial z_i\partial \bar{z}_j}\Big)\,.
$$
One can show that the signature of $L$ is a biholomorphic invariant.
If for example $L$ is positve-definite, then supposing that $p=0$ 
and that the CR-tangent space is given by $z_1=0$ one can 
introduce holomorphic coordinates $(z_1,z')$ so that the restriction
of $\rho $ to $\{z_1=0\}$ is 
$$
\rho (z)=\Vert z'\Vert^2+O(3)\,.
$$
Thus locally near $p$ the CR-tangent space lies outside $D$ and
except at $p$ is contained in the complement of the closure of $D$.
One can think of it as a (local) supporting complex hypersurface outside
$D$ at the point $p$.  If the Levi-form is positive-definite at 
every point of $\partial D$, one says that the domain (or its boundary)
is strongly pseudoconvex. If at each point $p\in \partial D$ the Levi-form
is only positive semidefinite, one just says that $D$ is pseudoconvex.
One can imagine that there is a huge difference between the concepts, 
particularly if the rank of $L_p$ is allowed to vary wildly
with the point $p$.  

\bigskip\noindent
In the early part of the 20th century E. E. Levi realized that
if $L_p$ is not positive-semidefinite, then every function defined
and holomorphic on $D$ near $p$ extends holomorphically across 
$\partial D$ at $p$.  On the other hand, if it is positive-definite, 
locally in the coordinates used for the above normal form, the
function $\frac{1}{z_1}$ is holomorphic on $D$ near $p$ and
does not continue across the boundary.  Therefore one asks
if the same holds at the global level.  Levi himself showed
that if at some $p\in \partial D$ the Levi-form is not positive
semidefinite, then every function holomorphic on $D$ continues across 
$\partial D$.  The Levi-Problem can be stated as follows:

\medskip
\centerline{$D$ pseudoconvex $\overset{?}{\Rightarrow}$ $D$ is a domain
of holomorphy.}
\medskip

In other words, if $D$ is pseudoconvex, given a divergent sequence
$\{z_n\}$ one would like to prove that there exists a holomorphic function
$f$ on $D$ with $\lim \vert f(z_n)\vert =\infty $.  Recall that this
property is the only one of Stein's axioms that is not automatically
fulfilled for a domain in $\mathbb {C}^n$. Hence, for domains one
is really asking if pseudoconvex domains are Stein with the
above property being called \emph{holomorphic convexity}.

\bigskip\noindent
In 1942 Oka solved the problem for pseudoconvex domains 
unramified over $\mathbb {C}^2$, and in 1953/54 this was 
extended to arbitrary dimensions independently by Oka,
Bremermann and Norguet.  Grauert points out that by
using the Behnke-Stein theorem (limits of domains of holomorphy
are domains of holomorphy) in $\mathbb {C}^n$, one needs only
to prove the result for strongly pseudoconvex domains.  This
is not true for domains in arbitrary manifolds.  For example,
it is a simple matter to construct a domain $D$ in a torus
with Levi-flat boundary, i.e., the Levi-form of an appropriate
boundary defining function vanishes identically (nevertheless
$D$ is pseudoconvex!), with the property that every holomorphic
function on $D$ is identically constant.  Grauert has constructed
much more sophisticated examples in \cite{GB}, where except for
a small set $\partial D$ is strongly pseudoconvex. 

\bigskip\noindent
Grauert states his theorem for bounded domains $D$ with smooth
boundaries in arbitrary complex manifolds:

\medskip
\centerline{strongly pseudoconvex $\Rightarrow $ holomorphically convex.}
\medskip

Actually he proves much more: \emph{If $D$ is strongly pseudoconvex,
then it contains finitely many pairwise disjoint maximal compact
analytic subvareties which can be blown
down so that the resulting complex space is Stein.} For details see
our discussion of his paper \emph{\"Uber Modifikationen und
exzeptionelle analytische Mengen}.

\bigskip\noindent
Grauert's elegant proof begins by implementing 
his now famous bumping technique where he
constructs a (finite) increasing sequence $\{D_k\}$ of domains 
containing $D$ such that at each step the restriction mapping
$H^\nu(D_{j+1},\mathcal {O})\to H^\nu(D_j,\mathcal {O})$ is surjective
for all $\nu\ge 1$.  The largest domain contains $D$ as a relatively
compact subset. Hence he proves that for $D'$ sufficiently near
$D$ with $D\subset \subset D'$ the same surjectivity result holds.
Actually in the final step of his proof he uses this for the
sheaf of a holomorphic line bundle where the surjectivity is proved
in exactly the same way.
\bigskip\noindent
Having achieved the above indicated surjectivity Grauert applies 
L. Schwartz' Fredholm theorem which states
that if $\varphi :E\to F$ is a surjective, continuous linear map of 
Fr\'echet spaces and $\psi:E\to F$ is compact, then 
$\varphi + \psi $ has closed image of finite codimension.  To apply this
Grauert organizes a (finite) Leray covering $\mathcal {U}$ of 
$D'$ which is refined
to a Leray covering $\mathcal {V}$ of $D$ and considers the
mapping
$$
\varphi: Z^q(\mathcal {U},\mathcal {O})\oplus C^{q-1}(\mathcal {V},\mathcal {O})
\to Z^q(\mathcal {V},\mathcal {O})
$$
which is the sum of the restriction map $R$ and the $\check{\rm C}$ech 
boundary map
$\delta$.  The cohomological surjectivity implies that this map is surjective.
Hence, the Schwartz Theorem implies that the image of $\delta =\varphi-R$ has
finite codimension. This is exactly the desired finiteness theorem.

\bigskip\noindent
Grauert uses this finite dimensionality to prove a result that is actually
stronger than the holomorphic convexity of $D$: Given a boundary point
$x_0$, he enlarges $D$ as above so that in addition $D'$ contains
a $1$-codimensional complex submanifold $S$ which contains $x_0$ but
is entirely contained in the complement of $D$. Using the finite
dimensionality of the cohomology of powers $F^k$ of the line
bundle defined by $S$ as well as its restriction to $S$, for $k$
sufficiently large he finds a section $s$ of $F^k$ which does not
vanish at $x_0$.  Thus if $h$ is the defining section of $S$
in $F$, then $\frac{s}{h^k}$ is a meromorphic function on $D'$
which is holomorphic on $D$ with a pole at $x_0$. 

\bigskip\noindent
As the title of the article indicates, Grauert applies his theorem
to prove that paracompact real analytic manifolds $M$ can be
embedded in Euclidean spaces of the expected dimension. For
$M$ compact this result was proved by Morrey using PDE methods
slightly earlier.  At the time Bruhat and Whitney had shown
that $M$ can be regarded as the set of real points of a complex manifold
$X$, and Grauert then constructed a non-negative strictly 
plurisubharmonic function
$\rho $ on a neighborhood of $M$ in $X$ which vanishes exactly on $M$ and
otherwise has non-vanishing differential. If $M$ is compact, then
Grauert's solution to the Levi problem immediately implies that
$T=\{\rho <\varepsilon\}$ is Stein and it is immediate that there
is an everywhere maximal rank, injective holomorphic map of $T$
(shrunk a bit) to some $\mathbb {C}^N$.  If $M$ is not compact,
more sophisticated arguments must be used. In particular, in order
to prove that there is a Stein tube one must use a generalized version
of the Behnke-Stein Runge theorem which is proved in \cite{GD} (see
our discussion of \emph{Levisches Problem und Rungescher Satz}). 
Then Remmert's theorem for Stein manifolds yields the
desired embedding.   

\bigskip\noindent
In closing we should add that for weakly pseudoconvex domains in
Stein manifolds, the Docquier-Grauert results are optimal.  Strongly
pseudoconvex domains in Stein spaces have been handled by Narasimhan 
\cite{NL} with the analogous results to those discussed above. 
There are numerous partial results for weakly pseudoconvex domains,
also in the case where singularities play a role. However, the general
situation is far from being understood.

\bigskip\noindent
\textbf{Levisches Problem und Rungescher Satz
f\"ur Teilgebiete Steinscher Mannifaltigkeiten}

\medskip\noindent
Here Grauert and Docquier begin by discussing nine conditions 
which are relevant for the study of the pseudoconvexity of a 
complex manifold. These are denoted by ($h$, $p_1,\ldots , p_7, p_7^*$)
and are interrelated by a graph
of implications with the condition $h$ of holomorphic convexity
being the strongest and $p_7^*$ being the weakest.  The latter
condition is a weak version of the condition that a Hartogs figure
cannot be mapped biholomorphically into the manifold so that
its image is not relatively compact but the image of its Shilov
boundary is compact. Riemann domains $G$ which are
unramified over a Stein manifold $M$ are considered,
and it is shown that if such a domain satisfies $p_7^*$, then
it is Stein. It is therefore holomorphically convex, i.e., $h$
is satisfied and consequently all of the conditions are fulfilled.

\bigskip\noindent
The basic idea of the proof is to apply Remmert's embedding theorem
to embed $M$ in some $\mathbb {C}^n$.  Then, using
the normal bundle of $M$ in this embedding the domain
$G$ is thickened to an unramified domain $\hat {G}$ over $\mathbb {C}^n$ which
also satisfies $p_7^*$. In that situation Oka's methods can
be applied to $\hat {G}$ to achieve the desired result.

\bigskip\noindent
Again using the idea of thickening a Remmert embedding, Runge approximation
theorems are proved via Oka-Weil approximation for domains unramified
over $\mathbb {C}^n$.  A condition for a Stein domain $M$ to be
Runge in a complex manifold $\check {M}$ (strongly simplified
for our presentation) is that there is a continuous increasing family of
Stein domains $M_t$ starting at $M$ and ending 
at $\check{M}$.  For $t_1<t_2$ it is proved that
$M_{t_1}$ is Runge in $M_{t_2}$ and then by
the classical Runge Theorem of Behnke and Stein it follows that
$\check {M}$ is Stein.  

\bigskip\noindent
This last mentioned result is really just a corollary of results 
proved in much greater generality.  However, we wanted to particularly
underline it, because it is exactly what is needed in \cite{GL} for proving
that Grauert's tube around a \emph{non-compact} real analytic
manifold is Stein.

\bigskip\noindent
\textbf{Th\'{e}or\`{e}mes de finitude pour la cohomologie des
espaces complexes}

\medskip\noindent
Recall that a smooth function on, e.g., a domain in
$\mathbb {C}^n$ is strictly plurisubharmonic if its
Levi form (complex Hessian) is positive-definite.
Stein manifolds are those complex manifolds $X$ which
possess a strictly plurisubharmonic exhaustion. This
is, in a certain sense, the solution of the Levi problem.  
In the work of Andreotti-Grauert, which we will now review,
a $q$-Levi problem is solved.  Just as the Levi problem
for strongly pseudoconvex domains was handled by Narasimhan
in the case of complex spaces, the context here is also for
complex spaces. However, in order to explain the essential
ideas, it is enough to consider the smooth case.

\bigskip\noindent
Let $B$ be a bounded domain with smooth boundary $\partial B$
in a complex manifold $X$. Andreotti and Grauert say that
$\partial B$ (or $B$) is $q$-pseudoconvex if the Levi form
of a defining function has at least $n-q+1$ positive eigenvalues,
the case of $q=0$ being reserved for compact manifolds.  The 
main goal of the paper is to prove the finite dimensionality  
of $H^k(B,\mathcal {F})$ for $k\ge q$ where $\mathcal {F}$ is 
a coherent sheaf on $X$.  For example, if $\rho :X\to \mathbb {R}^{\ge 0}$
is an exhaustion which is $q$-pseudoconvex outside of a 
compact set $K$ which is contained in a $\rho $-sublevel set
$B$, then the restriction map $H^k(X,\mathcal {F})\to H^k(B,\mathcal {F})$
is an isomorphism.  Thus the finite dimensionality for $X$
follows and if $K$ is empty, then $H^k(X,\mathcal {F})=0$ for
$k\ge q$.  

\bigskip\noindent
A domain $B$ is strictly $q$-pseudoconcave if $\rho $ is a defining
function as above and $-\rho $ is strictly $q$-pseudoconvex.
With the addition of some technicalities in the case of singular
spaces and coherent sheaves which are not locally free, the 
finiteness and vanishing theorems hold in the concave case for 
$0\le k\le n-q$.
  
\bigskip\noindent
It is interesting that, in order to prove these global results,
the main new work needed is of a local nature!
This point actually comes up in Grauert's previous paper, but since
it is handled by classical methods, one tends to forget it.
In that paper $\partial B$ is strongly pseudoconvex. 
In a Stein coordinate chart $U$ containing a given boundary point,
Grauert constructs a bump on $B$ to obtain a domain $B_1$.  It
is of fundamental importance that the cohomology of $U\cap B_1$
vanishes, i.e., Cartan's Theorem B for this domain!  Here, in
the $q$-pseudoconvex case, the analogous ``Lemma'', along with
Grauert's proof idea of bumping and then using the 
Fredholm Theorem of L. Schwartz, yields the proof.  Of course
there are substantial technical preparations which must be carried out.

\bigskip\noindent
Roughly speaking the above mentioned Lemma amounts to doing the
following. In a local coordinate chart $U$ at a boundary point
$\xi_0$ one chooses a transversal polydisk of dimension $n-q+1$
so that the restriction to it of the boundary defining function
is strictly plurisubharmonic.  Then one thickens it to obtain
a holomorphic family of such polydisks parameterized by, e.g.,
a polydisk of complementary dimension. Then, as in the
strongly pseudoconvex case, one creates a bumped region $B_1$ 
which intersects each transversal polydisk in a strongly pseudoconvex 
region and which only changes $B$ in a compact region in
$U$.  It is now necessary to prove a cohomology vanishing theorem 
for $V:=U\cap B$ for $k\ge q$, and the authors do exactly this by
viewing $V$ as a family of $(n-q+1)$-dimensional Stein domains. 
One of the main difficulties for this is proving
the appropriate Runge theorem.

\bigskip\noindent
\textbf{Algebraische K\"orper von automorphen Funktionen}

\medskip\noindent
Although the discussion here is in fact quite general,
applying to any pseudoconcave complex space, the work
in this paper is carried out in the special situation
where modular and associated automorphic functions are
playing an essential role.  In fact only one example
is considered, the quotient of the Siegel upper halfplane
$H$ by the modular group $\Gamma $, but as Borel later
pointed out, there is a wide class of examples where
the Andreotti-Grauert method applies. 

\bigskip\noindent
A (connected) complex space $X$ is said to be pseudoconcave if it 
contains a relatively compact open subset $Y$ with the property
that for every $x\in \partial Y$ there is a map
$\varphi :\mathrm {cl}(\Delta) \to \mathrm {cl}(Y)$ which is
holomorphic in a neighborhood of the closure of a $1$-dimensional
disk with image in the closure of $Y$ with the properties that
$\varphi (0)=x$ and $\varphi (\partial \Delta)\subset Y$. By 
thickening such ``disks'' Andreotti and Grauert obtain a double
covering of $\mathrm{cl}(Y)$ by images of polydisks (one relatively
compact in the other) so that their Shilov boundaries are contained
in $Y$.  Here the Shilov boundary of a polydisk 
$\Delta =\{\vert z_i\vert < 1,\, i=1,\ldots ,n\}$ is the
set where $\vert z_i\vert =1$ for all $i$.  They then apply Siegel's
method using the classical Schwarz Lemma to prove the following
fact:  The field $\mathbb {C}(X)$ of meromorphic functions on $X$
is a finite algebraic extension $\mathbb {C}(f_1,\ldots ,f_k)[g]$
of the field of rational functions in $k$-algebraically independent
meromorphic functions where $k\le \mathrm{dim}(X)$. It should be
mentioned that Andreotti went on to develop this theory in 
several ensuing works.

\bigskip\noindent
If, for example, $X$ arises as the quotient $\widehat {X}/\Gamma $
of some other space by the proper action of a discrete group,
then the notion of pseudoconcavity can be formulated at
the level of $\widehat {X}$.  Andreotti and Grauert do this
and then restrict their attention to the case where 
$\widehat {X}=H$ is the Siegel upper halfplane of complex
$n\times n$-matrices $Z=X+iY$ which are symmetric and where
$Y>0$.  The discrete group which is of interest here
is $\Gamma =\mathrm{Sp}_{2n}(\mathbb {Z})$. It is acting properly
and discontinuously so that the quotient $H/\Gamma $ has the
natural structure of a complex space. It is well known that 
$\Gamma $-periodic meromorphic functions, i.e., functions on the
quotient, are important in more than one area of mathematics.

\bigskip\noindent
Using a well-known fundamental region $\Omega_0$ for the $\Gamma $-action
along with the strictly plurisubharmonic function 
$k(z)=-\mathrm{log}\vert Y\vert$, Andreotti and Grauert determine
a region in $H$ that descends to the quotient to show that it 
is pseudoconcave. An essential part of the proof is devoted to achieving
the periodicity of $k(z)$ by minimizing it over $\Gamma $.  This
can be done, because the minima are taken on in the fundamental region.

\bigskip\noindent
As a consequence, the result on function fields can be applied in this
case. This was known already, but was proved by using vastly more complicated
methods. Furthermore the possibility of using pseudoconcavity in 
this area of mathematics was a totally new, extremely useful 
idea. It should be remarked that the same type of method can be
used to prove that natural spaces of modular forms, e.g., for
the canonical bundle, are finite dimensional. Furthermore, pseudoconcavity
implies that the quotients $H/\Gamma $ close up in projective embeddings
to compact complex spaces to which all meromorphic functions extend.

\bigskip\noindent
\textbf{\"Uber Modifikationen und exzeptionelle analytische Mengen}

\medskip\noindent
Given a complex space $X$ and a compact subvariety $A$ one
is interested in understanding when there is a complex
space $Y$ with a distinguished point $y\in Y$ and a surjective
holomorphic mapping $\pi :X\to Y$ which is biholomorphic
from $X\setminus A$ to $Y\setminus \{y\}$ and with $\pi(A)=\{y\}$. 
In other words, one would like to have sufficient conditions
for $A$ to be blown down to a point.  In the projective algebraic
setting in the case where $X$ is a surface 
certain results were already known, e.g., for blowing down
a smooth rational curve.  In this beautiful paper Grauert answers this
question in a general analytic setting in terms of the neighborhood
geometry of $A$ and its normal bundle. Underway he proves a number
of results that can be considered as preparatory but which are also
extremely useful in many areas of global complex geometry. As is often
the case for Grauert, the guiding light is 
given by the notion of strong pseudoconvexity.  

\bigskip\noindent
Here Grauert begins by noting that his solution to the Levi problem
for relatively compact domains $G$ in complex manifolds $X$
had just been extended to the case where $X$ is singular
(\cite{NL}).
Using Remmert's reduction theorem, he observes that the result
can be stated as follows: If $G$ has strictly pseudoconvex
boundary, then it contains a maximal compact analytic subset
$A$ which can be blown down to a finite number of points 
(corresponding to its connected components) by a map $\pi:G\to Y$
where $Y$ is a Stein space.  Conversely, if a connected
compact analytic set can be blown down to a point, then it
has a strongly pseudoconvex neighborhood.  So it is natural
to study the relation of this type of question to the
pseudoconvexity of neighborhoods of the $0$-section of the
normal bundle of $A$ or more generally for any bundle.  

\bigskip\noindent
For line bundles $F$ over a compact complex manifold $X$ the importance
of the notion of the positivity of a Hermitian bundle metric
was known. One says that $F$ is ample if some power $F^k$ defines an embedding
$X\hookrightarrow \mathbb {P}(\Gamma (X,F^k)^*)$ by mapping a
point $x\in X$ to the hyperplane of sections which vanish at $x$.
Kodaira's basic theorem states that $F$ is ample if and only if 
it possesses a positive bundle metric.  Since the region defined 
by $\Vert \cdot \Vert >1$
can be regarded as a tubular neighborhood of the $0$-section
of $F^*$, Grauert reformulates positivity in terms of the
strong pseudoconvexity of the $0$-section of the dual bundle.
It is important to emphasize that this also makes sense in the
case where $X$ is singular.
He calls this property \emph{schwach negativ}. Thus 
the embedding theorem can be stated as $F$ is ample if and only
if $F^*$ is schwach negativ. This is then equivalent to the 
$0$-section of $F^*$ being the maximal compact subset of $F^*$. It 
can be blown down to a Stein space which Grauert shows to be
affine.  Grauert's notion for vector bundles 
of higher rank is defined analogously: A vector bundle $V$ over
a compact complex space is said to be Grauert-positive if and 
only if the dual bundle $V^*$ is schwach negativ in the above sense,
i.e., its $0$-section can be blown down. It should be remarked that 
in the vector bundle case the relation of Grauert's positivity condition to 
Griffiths-positivity is still not understood.

\bigskip\noindent
One of the main results of the paper is the embedding theorem:
If a complex compact space $X$ possesses vector bundle which
is schwach negativ, then it is projective algebraic.  The key is the
Stein property for the blown down bundle space. Even in the
case of line bundles $F\to X$ the result is new, because here singular spaces
are allowed.  This also gives a proof of the embedding theorem
for Hodge spaces, another result of Kodaira in the smooth case.

\bigskip\noindent
The following is a more general version of the fact 
mentioned above, i.e., that the
blown down dual bundle space is affine:  Let
$F$ be the bundle of a divisor which has support $A$.  Suppose
$F\vert A$ is positive and that $X\setminus A$ contains no
positive dimensional compact analytic subsets.  Then $X\setminus A$
is affine and $F$ is a positive bundle on $X$.  The key ingredient 
for the proof, called a Hilfslemma by Grauert, is probably even more
useful: A line bundle $F$ over a compact complex space $X$ is
positive if and only if for every analytic subset $A$ there exists
$k>0$ so that $F^k\vert A$ has a section which vanishes at some point
of $A$ but does not vanish identically. 

\bigskip\noindent
Returning to the main theme of the paper, Grauert considers the notion
of the normal bundle of a compact complex subvariety $A$ of a complex
space $X$.  Due to the possible singular nature of these spaces,
this must initially be regarded as a sheaf corresponding to the 
ideal sheaf $\mathfrak {m}$ of $A$ or more generally any coherent 
ideal sheaf $\mathcal{I}$ which defines $A$.
In typical Grauert fashion 
he is not phased by this difficulty
but rather introduces the (quite natural) notion of a linear fiber 
space associated to a coherent sheaf.  Locally over a trivializing
neighborhood $U$ this is a subvariety of $U\times \mathbb {C}^n$
where the fibers are subvector spaces of $\mathbb {C}^n$ so that
addition and scalar multiplication are well defined. Thus, given
$A$ and the ideal sheaf $\mathcal {I}$ as above one has its normal 
linear fiber space $N_{\mathcal {I}}$
with its $0$-section and the notion of schwach negativ has the
obvious meaning. In elegant fashion Grauert transfers the pseudoconvexity
of a neighborhood of the $0$-section to that of a neighborhood of
$A$ in $X$ and proves the desired result: $A$ can be blown down
if for suitable $\mathcal {I}$ the normal linear fiber
space $N_{\mathcal {I}}$ is schwach negativ.

\bigskip\noindent
Of course the results in this paper have numerous applications.
Even in the case of surfaces one needs Grauert's
results to show that an irreducible curve $C$ 
has negative self-intersection number if and only if it can 
be blown down to a point.  Grauert's most general result in this direction
is that a $1$-dimensional subvariety in a surface can be blown down
if and only if its self-intersection matrix is negative definite.
\subsubsection* {Direct image theorem}
\textbf{Ein Theorem der analytischen Garbentheorie und die 
Modulr\"aume komplexer Strukturen}

\medskip\noindent
The proof of the Direct Image Theorem (\emph{Bildgarbensatz}) (\cite{GDI}) is
one of Hans Grauert's greatest accomplishments.  We will state it
here, say a bit about the proof and give an application mentioned
by Gauert in the paper.  As we wrote in (\cite{H1}), the applications
are so far reaching in complex analytic geometry that it would be 
unimaginable to work in the area without having it available.  

\bigskip\noindent
Let us turn to the setting of complex analysis at the time (the late 1950's).
A great deal was known, at least compared to ten years before.  The
notion of a complex space had been clarified, Grauert already had a
huge experience as \emph{Handarbeiter} in dealing with problems of
cohomology, e.g., the Oka Theorems and his proof of Theorems A and B
were behind him, and he understood very well how to deal with refining
covers and using the relevant Fr\'echet spaces and compact 
operators between them.  Given all of this he was in a position
to consider the problem of the coherence of direct images of
coherent sheaves.

\bigskip\noindent
The initial geometric context of this theorem is quite simple. One
begins with a holomorphic map $F:X\to Y$ between complex spaces
which is defined in the most naive way, e.g., locally it is
given by holomorphic functions.  Associated to an open subset $U$
in $Y$ one has the algebra $\mathcal {O}_X(F^{-1}(U))$ on its
preimage. This has the structure of an $\mathcal {O}_Y(U)$-module
which is given by multiplication by lifted functions $F^*(f)$.  This
presheaf defines a sheaf $F_*(\mathcal {O}_X)$ on $Y$. It contains
a great deal of information about the map $F$, in particular about
its singularities. It would clearly be of interest to know whether
or not it is coherent.

\bigskip\noindent
If indeed the direct image $F_*(\mathcal {O}_X)$ is coherent, then
its support is a closed analytic subset of $Y$.  Hence, the correct
condition for the direct image theorem to hold must be something
that guarantees that images of analytic subsets are analytic.
At the time, Remmert's theorem, which guarantees that this is
the case for $F$ being a proper (holomorphic) map, had been proved. 
It should be underlined that even the notion \emph{proper}, i.e., inverse
images of compact sets are compact, was rather new.  Cartan had
introduced this in the 1930's while discussing the fact that
the action of group of automorphisms on a bounded domain is
proper and the notion was explicitly described in Bourbaki.
Kuhlmann had pointed out that there is a weaker notion (semi-proper)
and had proved that Remmert's theorem holds for this kind of map.
Stein and Grauert were always interested in understanding holomorphic
equivalence relations and finding a good condition which
would insure that the quotient is analytic. In fact, one of Grauert's 
last papers was devoted to a situation where a sort of semi-properness
was built into the assumptions.

\bigskip\noindent
In any case, at the time when Grauert considered the problem of the
coherence  of direct images of coherent sheaves, much was known,
but even at the set-theoretic level (Remmert's theorem) things
had not settled in. There were also a huge number of foundational
issues.  For one, even the notion of an analytic morphism had to
be improved. One reason for this, at least from Grauert's point of view,
was that the entire project had to be carried out in the context of
complex spaces where the structure sheaf is allowed to have \emph{nilpotent}
elements.  This means that the local model is as before an analytic set $A$
in some domain $D$ in $\mathbb {C}^n$, but the sheaf of germs of holomorphic
functions is $\mathcal {O}_D/\mathcal {I}_A$ where $\mathcal {I}_A$ is
any coherent ideal sheaf which defines $A$ as its $0$-set.  Since the
structure sheaf is not necessarily a subsheaf of the sheaf of continuous
functions, the classical definition of a map being holomorphic, i.e.,
pullbacks of holomorphic germs are required to be holomorphic, is not
sufficient.  A holomorphic map is then a pair $(F_0,F_1)$ where
$F_0:X\to Y$ is a usual map of sets and $F_1$ is a map of structure
that encodes the notion of pullback, a continuous homomorphism of
sheaves of algebras  $F_1:Y\times_{F_0}\mathcal {O}_Y\to \mathcal {O}_X$.  
Grauert begins his paper with a rather long discourse on how to deal with these
\emph{new} complex spaces where he proved the key theorems such
as Theorems A and B in this more general setting.

\bigskip\noindent
Given a morphism $F:X\to Y$ and sheaf $\mathcal S$ of 
$\mathcal {O}_X$-modules, $F_1$ is applied to equip the
direct image sheaf $\pi_*(\mathcal {S})$ with the structure
of a sheaf of $\mathcal {O}_Y$-modules.  One can go an important
step further:  For  $U$ open in $Y$ and every $q\ge 0$ the cohomology space 
$H^q(F^{-1}(U),\mathcal {S})$ is equipped by means of $F_1$ with the
structure of a $\mathcal {O}_X(U)$-module.  Hence, for every $q$ we 
have direct image sheaf $R^qF_*(\mathcal {S})$.  The following is then
the \emph{Bildgarbensatz}.
\begin {theorem}
If $F:X\to Y$ is a proper holomorphic map of complex spaces and
$\mathcal {S}$ is a coherent sheaf of $\mathcal {O}_X$-modules,
then for every $q\ge 0$ the direct image sheaf $R^qF_*(\mathcal {S})$
is coherent.
\end {theorem}
There were germs of this result around at the time, e.g., Remmert
had proved a result in special case of finite maps, 
Grauert and Remmert had proved it
in the situation where $X=Y\times \mathbb {P}_n$ and $F$ is the
obvious projection and Grothendieck had proved the far simpler algebraic
version.  However, this result, along with the Oka Principle papers,
brought complex analysis into a new era.  Only five years before members
of the Behnke seminar were trying to understand $\sqrt{xy}\,$!

\bigskip\noindent
Let us quote Grauert (\cite{GSW}, p. 446) when discussing the main
difficulties in the proof which involve a power series argument
to handle the directions transversal to a fixed fiber $X_0$.
\emph{Roughly speaking, the proof of the direct image theorem uses power
series expansion whose coefficients are cohomology classes on one fixed
fiber. The coefficients are obtained recurrently and with estimates.
Cohomology spaces carry natural Fr\'echet space structures.  However, to
get a convergent power series by recurrent formula with estimates, one
needs a fixed norm for the iteration process instead of the infinite
sequence of semi-norms.  One key point of the proof of the direct image
is to replace a cocycle with estimates for a weak norm by another
cocycle with estimate for a stronger norm modulo a coboundary with estimates
for an even weaker norm.}

\bigskip\noindent

Grauert began this article with a rather lengthy discussion of
the importance of the cohomology of a certain 
direct image sheaf for the study of moduli spaces $M$.
In that case $F:X\to M$ is a usual (surjective) proper holomorphic map 
of complex manifolds which is everywhere of maximal rank.  The 
relevant sheaf is $\Theta_X$, the sheaf of germs of holomorphic vector
fields on $X$. For example, if $X_y$ denotes the fiber over $y\in M$
and $r_q(y)$ is the dimension of $H^q(X_y,\Theta_y)$, then the 
semicontinuity of $r_q(y)$, which was proved by Kodaira and Spencer
using the method of harmonic integrals, follows immediately from
the direct image theorem.  In fact, in much greater generality 
the direct image theorem implies
the semicontinuity for any coherent sheaf provided the proper
map $F:X\to Y$ is flat.

\bigskip\noindent
It took the complex analysis community a number of years to understand
and somewhat simplify Grauert's proof (see, e.g.,\cite{FK} and \cite{N}). 
According to Grauert, the simplest proof can now be found in \cite{GR2}.

\subsection* {Akademischer Lehrer}
Let me close this note with some personal comments. I was introduced
to the ``German school'' of complex analysis in the late 1960's
in the Stanford lectures of Aldo Andreotti, where I was the only
student.  In the first semester of these lectures Andreotti explained
a number of Grauert's results (some of them discussed above) in a
beautiful way.  Although I had been a student for a while, this was
the first time that I felt that I had seen ``the truth''.  Of course
Andreotti was a master lecturer, but the truth, I sensed, was embedded
in Grauert's work.  A few months later, when I realized I should
prepare for my German language exam, luck struck again:  My advisor,
Halsey Royden, had explained in seminars some of his ideas on
metrics on Teichm\"uller space. As a result 
I optimistically thought it might be good
to go back to the basics and read Hermann Weyl's book \emph{Die
Idee der Riemannschen Fl\"ache} and then jump to the modern
developments and study Grauert's paper \emph{\"Uber Modifikationen
und exzeptionelle analytische Mengen}. Up to this point I had had very
little experience reading mathematics and assumed this is the way
it should be! Looking back, it is hard for me to believe that I had
been so naively audacious!  A year or so later 
during my postdoc time in Pisa, I told
my friends that these two works were the only things I really
understood.  At the time I was embarrassed to admit this  
but soon realized my good fortune!

\bigskip\noindent
In the second semester of the above-mentioned course Andreotti ``asked''
me to lecture on various topics involving $\bar \partial$ at the 
boundary, Hans Levy's extension, etc.. The audience consisted of 
Andreotti and Wilhelm Stoll.  This began
my lasting friendship with Stoll.  He came from another ``Schwerpunkt''
of complex analysis, namely from the T\"ubingen group of Hellmuth Kneser.
The complex analysis of T\"ubingen was, to a certain extent, related
to that of the M\"unster school, e.g., they had competing theories
of meromorphic maps.  However, Kneser and Stoll went in other directions,
proving continuation theorems under assumptions of bounded volume,
and then building the foundations of value distribution theory
in several complex variables.

\bigskip\noindent
My close relationship with Stoll continued during the almost 
10 years I spent at Notre Dame where, particularly due to
Stoll's connections, the faculty had close ties to 
German mathematicians.
In my very first year there Stein visited for a semester. His 
energy, openness and obvious love of mathematics made a great 
impression on me.  The next year Remmert came,
and it was a great honor for me to drive him around town as he reviewed
periods in which he had also been a guest at Notre Dame.  
In those days it was a bit non-standard to go from the US to
Oberwolfach for a week, but when he invited me I jumped at the
opportunity.   In my Oberwolfach lecture Grauert, Remmert, Stein and Forster
were in the first row. I figured if I could get through that
I could get through anything!

\bigskip\noindent
Grauert was a kind, warm person of very few words. Sometimes he looked
formal, but he was not. In the above-mentioned conference Douady lectured 
on his construction of the versal deformation of a complex space.  
Grauert, who had developed his own version of this theory,
sat in the first row.  Douady, who was dressed in a silk-like Hawaian shirt
which was not completely buttoned and was rather dirty because
the night before he had slept in the forest, explained his puzzles and
made silly jokes. Grauert remained quiet and respectful, only
caring about the content. Later on at a memorial meeting in honor of 
Andreotti, who had passed away at a very early age, 
one could see that Grauert and Douady were very close.

\bigskip\noindent
Grauert didn't say much, but when he did he meant it!
A comment of ``good, continue on'' to a young speaker after
a talk really meant something. Similarly, his way of praising a student
was often ``Das k\"onnen wir so machen''. He was a key referee
for one of our research concentrations sponsored by the DFG. 
As a 50-year-old I nervously appeared in front of Grauert 
for his comments: ``almost everything
is good, but the mathematics in this subproject is not important''. 
Of course I dropped the subproject from the proposal.
One might think of Grauert as being opinionated, but his opinions
were based on serious thought; anything he said or wrote should
be taken seriously!  

\bigskip\noindent
Hans Grauert was an ``Akademischer Lehrer'' in the sense of Humboldt.
He didn't teach a subject because it was in the syllabus; he taught
it because he had thought about it and knew it was important. He
had a large number of doctoral students, more than 40, and he was
proud of it. I know that in every case he had thought through
their projects and made notes on the little pieces of paper that
he carried around.  One of our colleagues who was a student of Grauert
remembers seeing the same piece of paper every time he came to
Grauert's office hours.

\bigskip\noindent
Even though he was not the type of person to be heavily 
involved with the global organization of science, he did his duty, 
e.g., as managing editor of Mathematische Annalen and working in 
various capacities with the DFG. He certainly respected the 
traditions of G\"ottingen
and was proud to have been President of the G\"ottingen Academy 
of Sciences.  His strong points were, however,
in the classroom where he focused on important phenomena in
mathematics, in his one-on-one work with his students and of course
in his remarkable research.  

\bigskip\noindent
Those of us who have had the privilege of knowing Hans Grauert will not 
forget him. Fortunately, his deep ideas have survived him in his
written works.  Let us hope that his high standards
of excellence in every aspect of our science will be carried on by
future generations.   
\bigskip\noindent
Alan Huckleberry\\
Ruhr Universit\"at Bochum\\ 
and\\
Jacobs University in Bremen\\
ahuck@cplx.rub.de
\end {document}